\documentclass[letterpage, 11pt, notitlepage]{article}
\usepackage[margin=1.0in]{geometry}

\usepackage{hyperref} 

\usepackage{times,tabularx,subfigure}
\usepackage{amssymb,amsfonts,amsmath,graphicx,cite,caption}
\usepackage{enumerate}
\usepackage[shortlabels]{enumitem}
\usepackage{multicol}
\usepackage{empheq}
\usepackage{color}

\newcommand{\calB}{\mathcal{B}}
\newcommand{\calG}{\mathcal{G}}
\newcommand{\calL}{\mathcal{L}}

\newcommand{\vct}[1]{\boldsymbol{#1}}

\newtheorem{thm}{Theorem}

\title{QCQP-Net: Reliably Learning Feasible Alternating Current Optimal Power Flow Solutions Under Constraints}

\author{
Sihan Zeng\thanks{J.P. Morgan AI Research}
\and Youngdae Kim\thanks{ExxonMobil Technology and Engineering Company -- Research}
\and Yuxuan Ren\thanks{Department of Computational Applied Mathematics \& Operations Research, Rice University}
\and Kibaek Kim\thanks{Argonne National Laboratory}
}

\begin{document}

\maketitle

\begin{abstract}%
At the heart of power system operations, alternating current optimal power flow (ACOPF) studies the generation of electric power in the most economical way under network-wide load requirement, and can be formulated as a highly structured non-convex quadratically constrained quadratic program (QCQP). Optimization-based solutions to ACOPF (such as ADMM or interior-point method), as the classic approach, require large amount of computation and cannot meet the need to repeatedly solve the problem as load requirement frequently changes. On the other hand, learning-based methods that directly predict the ACOPF solution given the load input incur little computational cost but often generates infeasible solutions (i.e. violate the constraints of ACOPF). In this work, we combine the best of both worlds -- we propose an innovated framework for learning ACOPF, where the input load is mapped to the ACOPF solution through a neural network in a computationally efficient and reliable manner.
Key to our innovation is a specific-purpose ``activation function'' defined implicitly by a QCQP and a novel loss, which enforce constraint satisfaction. We show through numerical simulations that our proposed method achieves superior feasibility rate and generation cost in situations where the existing learning-based approaches fail.\looseness=-1

\end{abstract}

\section{Introduction}








As one of the most important problems in modern power system operations, the study of alternating current optimal power flow (ACOPF) focuses on finding the most economical power generation scheme under network-wide load requirement and physical transmission constraints. Mathematically, ACOPF can be formulated as non-convex quadratically constrained quadratic program (QCQP) problem with the number of decision variables and constraints scaling proportionally with nodes and transmission lines in the power grid. The most common approach for reliably solving the ACOPF problem is limited to classical nonlinear optimization algorithms, such as interior-point method, 
which are highly computationally expensive for modern large-scale power systems involving at least thousands of nodes.
Furthermore, the constant fluctuations of the loads and conditions of the transmission line in addition to the uncertainty of energy supplies of renewable energy resources require the ACOPF problem to be solved repeatedly online, making the current optimization-based methods limited in real life.

With the recent advances in deep learning infrastructure and the improved ability to collect and store data, learning-based approaches have been proposed to solve complex optimization problems. 
The first attempt to solve ACOPF with machine learning has been made in \cite{guha2019machine}, where the authors employ a simple feed-forward neural network to parameterize the mapping from the input to the output of the ACOPF problem. However, the special structure of the ACOPF problem presents a peculiar challenge. The constraints define a nonlinear non-convex feasibility set around the optimal solution; while the neural network can consistently generate outputs close to the optimal solution in the Euclidean distance, they are not guaranteed to lie within the constraint set. In other words, learning-based approaches often produce highly infeasible solutions that cannot be directly deployed.

To address this issue, various methods have been proposed to encourage the output of the neural network to obey ACOPF constraints \cite{singh2021learning,pan2020deepopf,donti2021dc3}. Specifically, \cite{singh2021learning} enforces constraint satisfaction by adopting the Sobolev training scheme \cite{czarnecki2017sobolev}, which penalizes the mismatch in the Jacobian matrix at the solution in addition to the prediction error. 
Another line of works \cite{pan2020deepopf,donti2021dc3} recognizes that the solution of an ACOPF problem can be divided into 1) independent variables that control the power system operation and 2) dependent state variables that can be determined from control variables by solving the power flow equations. They propose predicting the control variables while leveraging additional loss functions to penalize constraint violation on the resulting state variables. The constraint violation loss is not readily differentiable, and different approaches such as zeroth order gradient or implicit function theorem are taken in these works to estimate/derive its gradient.

Although \cite{pan2020deepopf,donti2021dc3} significantly improves the constraint satisfaction, they build on the assumption that for any given control variable produced by the neural network they can always find a feasible solution (state variables) satisfying the power flow equations.
This assumption may not hold under high demand fluctuations, in which case the trained neural networks may be unreliable and cannot provide a meaningful solution.
To mitigate this issue, we solve a relaxation of the power flow equations with a focus on minimizing the constraint violations of these equations, which is also desirable for a reliable and robust operation of power systems. 

To achieve this, we formulate our relaxation problem as a non-convex QCQP and integrate the QCQP solver into the neural network as a differentiable activation function. 
We establish the differentiability of the non-convex QCQP activation function by extending the techniques first developed in \cite{amos2017optnet}, which only handles convex quadratic programs.
Equipped with properly designed loss functions, our proposed framework effectively deals with the infeasibility issue observed in current learning-based approaches such as \cite{pan2020deepopf} while showing similar competitive performance in feasible cases.

\subsection{Main Contributions}
In this work, our goal is to design an accurate and reliable end-to-end neural network architecture for predicting the solution of ACOPF while achieving high computational efficiency in both training and inference phases. 
%
As the first main contribution of our work, we propose a systematic pipeline and loss function for training a feed-forward neural network that maps the input load to the independent variables of the ACOPF solution. The state variables are then produced from the predicted independent variables by solving a relaxed variant of the power flow equations, which can be expressed as a non-convex QCQP. The relaxed power flow equations, as we discuss in details in Section~\ref{sec:learn_acopf}, are an important innovation of this work and allow us to train the neural network in a much more stable way when the control variables are imperfectly predicted.\looseness=-1



\begin{figure}[ht]
\centering
\includegraphics[width=.9\linewidth]{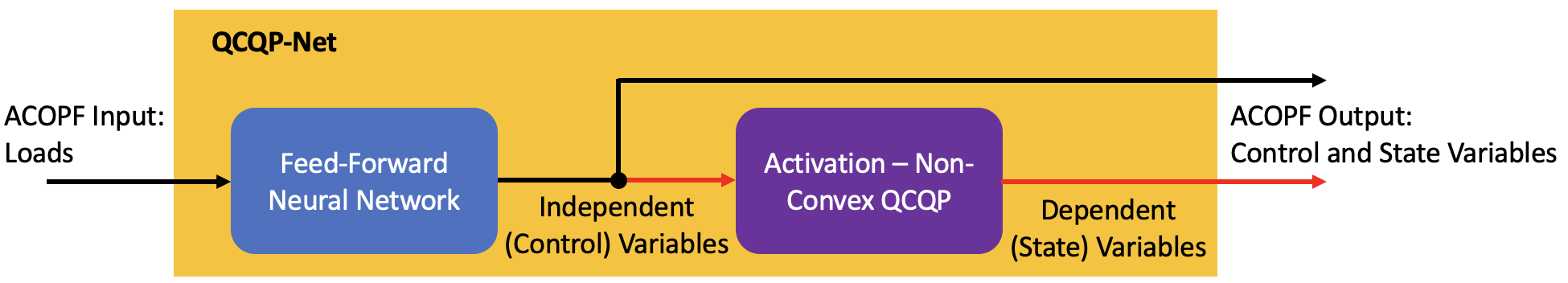}
\caption{QCQP-Net Architecture. Computation path in red only taken in training phase.}
\label{fig:QCQPNet}
\end{figure}

We can regard the relaxed power flow equations from a different angle as a specific-purpose activation function tailored to the ACOPF problem.
As a second main contribution, we establish the conditions under which this activation function is a differentiable mapping from the control variables to the state variables, and derive closed-form expressions of the (sub)gradient. This ensures that the downstream constraint violation loss on the state variables can be back-propagated. 
We name our proposed architecture QCQP-Net and show its structure in Figure~\ref{fig:QCQPNet}.
%
%
We numerically evaluate the performance of the proposed QCQP-Net on ACOPF problems of various scales in Section~\ref{sec:experiments}. The results show that in large power systems with wide load variations where the existing approaches fail to learn, QCQP-Net stably learns highly feasible solutions with low generation costs.

\subsection{Related Works}


This paper presents a novel learning framework specifically designed for reliably and efficiently solving ACOPF problems. It closely relates to the existing works that study ACOPF from both optimization and deep learning perspectives, and is inspired by recent advances in differentiable convex programming. 
We discuss the relevant literature in these domains to give context to our novelty.\looseness=-1

\noindent\textbf{ACOPF:} 
From an optimization perspective, a large volume of works seek to design provably convergent algorithms for ACOPF \cite{cvijic2012applications,kar2017admm,yuan2019second,sun2021two,wang2022nested} and to numerically accelerate classic nonlinear optimization solver through massively parallelized computation \cite{roberge2016optimal,huang2017performance,araujo2019simultaneous,kim2022accelerated,zhang2023solving}. In the learning regime, two main lines of work include 1) learning an end-to-end mapping from input of the ACOPF problem to the output \cite{singh2021learning,donti2021dc3,nellikkath2022physics,pan2020deepopf} and 2) learning parameters and/or sub-steps within an optimization solver \cite{baker2019learning,zhang2022learning,zeng2022reinforcement,sadat2021initializing}. While the former approaches allow for much faster inference, it may suffer higher constraint violation risk than the latter (which makes it less reliable and suitable for solving safety-critical power system problems). Our work is exactly motivated to address this issue.

\noindent\textbf{Differentiable Convex Programming:} Introduced in \cite{amos2017optnet} and later popularized by \cite{agrawal2019differentiable}, differentiable convex programming treats a convex optimization problem as an implicit definition of a mapping from the parameters of the optimization problem to the optimal solution. By carefully analyzing how a unit change in the parameters impacts the optimal solution through the lens of KKT conditions, \cite{amos2017optnet,agrawal2019differentiable} devise an innovative method for computing the (sub)gradients of the mapping. Applications of differentiable convex programming span neural network layer design \cite{amos2017optnet,agrawal2019differentiable,wang2019satnet}, inverse problems \cite{liang2019differentiable,mourya2023mcnet}, computer vision \cite{chen2020end,yeh2022total}, mechanism design \cite{curry2022learning,zeng2023near}, and many other domains. In this work we develop techniques to differentiable through a QCQP with quadratic equality constraints, which is much more challenging to handle due to the non-convexity. However, important pieces of our innovation are built upon \cite{amos2017optnet}.

\vspace{8pt}
\noindent\textbf{Outline of the paper.} The rest of the paper is structured as follows. In Section~\ref{sec:acopf}, we present the formulation of ACOPF and its important structure. In Section~\ref{sec:learn_acopf}, we propose a novel loss function for training an end-to-end prediction model for ACOPF by exploiting the problem structure. Evaluating the gradient of the loss function requires differentiating through a non-convex QCQP with quadratic equality constraint. We discuss how such differentiation can be performed in Section~\ref{sec:differentiable-qcqp}. Section~\ref{sec:experiments} presents the numerical simulations that demonstrate the stable and effective training of our end-to-end prediction model. 
Finally, we conclude in Section~\ref{sec:conclusion}.

\section{ACOPF Formulation}\label{sec:acopf}

We consider a power system represented by a connected graph with a set $(\mathcal{B},\mathcal{L})$ of buses and connected lines, respectively. Each node $i\in\mathcal{B}$, also referred to as a bus, has a complex power demand denoted as $d_i = p_i^d + j*q_i^d$ for some $p_i^d,q_i^d\in\mathbb{R}$. 
The voltage of bus $i$ is $v_i\in\mathbb{C}$, and we use $e_i$ and $f_i$ to denote the real and imaginary parts, i.e. $v_i=e_i+j*f_i$.
A subset of buses may have a power generator attached, and we use $\calB_\text{PV}\subseteq\mathcal{B}$ to denote the set of nodes with at least one generator attached. 
We use $\calG_i$ to denote the collections of generators attached to bus $i$ and define $\calG:=\cup_{i\in\calB} \calG_i$.
Each generator $g\in\calG$ can generate complex power with a real part $p_g\in\mathbb{R}$ and imaginary part $q_g\in\mathbb{R}$.


The edge of the graph, also referred to as a branch, represents a directed transmission line between two buses. For each branch $(i,j)\in\calL$ from bus $i$ to $j$, $p_{ij}$ and $q_{ij}$ denote the real and imaginary power flow in the normal direction. Power may also flow in the reverse direction, and we use $p_{ij}$ and $q_{ij}$ to denote the reverse power flow through branch $(i,j)\in\calL$. It is worth noting that $p_{ji}$ and $q_{ji}$ may not simply be the negative of $p_{ij}$ and $q_{ij}$ but are determined from the voltage at bus $i$ and $j$ by solving a system of power flow equations \eqref{eq:powerflow_1}-\eqref{eq:powerflow_4}, where parameters $B_{ij},G_{ij},B_{ji},G_{ji}\in\mathbb{R}$ are dictated by the physical properties of the power system.

For any bus $i$, we use $\mathcal{N}_i^{\text{from}}$ and $\mathcal{N}_i^{\text{to}}$  to denote its (directed) neighbors, i.e. $\mathcal{N}_i^{\text{from}}=\{j:(i,j)\in\calL\}$ and $\mathcal{N}_i^{\text{to}}=\{j:(j,i)\in\calL\}$.

The objective of the ACOPF problem, formulated in \eqref{eq:acopf}, is to find the most economic set points of generators that satisfy the power demand $p_i^d,q_i^d$ at every node $i$ under capacity limits and physical transmission laws. The generation cost function is quadratic in the real power output, where $c_{1,g},c_{2,g}\in\mathbb{R}_+$ are non-negative constant parameters for all $g\in\calG$. Eqs.~\eqref{eq:powerbalance_1}-\eqref{eq:powerbalance_2} are known as power balance equations and encode the power transmission laws along with Eqs.~\eqref{eq:powerflow_1}-\eqref{eq:powerflow_4}.
Eqs.~\eqref{eq:linelimits} state that the power flow magnitude between bus $i$ and $j$ cannot exceed the limit $\overline{s}_{ij}$. Eq.~\eqref{eq:voltagelimit} restricts the voltage at a bus to lie within a tolerable range. Eqs.~\eqref{eq:genlimits} represent the capacity of the power generators. The optimization problem can be expressed in a matrix form as a QCQP, but is obviously non-convex due to the quadratic equality constraints in Eqs.~\eqref{eq:powerbalance_1}-\eqref{eq:powerflow_4}.

\vspace{-5pt}
\begin{subequations}
\begin{align}
\min_{p_g,q_g,f_i,e_i,p_{ij},q_{ij},p_{ji},q_{ji}} \hspace{-30pt}&\hspace{35pt}
\sum_{g\in\mathcal{G}} ( c_{2,g} p_g^2 + c_{1,g} p_g ) \\ 
\text { s.t. }\quad 
& G_{ii} (e_i^2 + f_i^2) + \sum_{j \in \mathcal{N}_i^\text{fr}} p_{ij} + \sum_{j \in \mathcal{N}_i^\text{to}} p_{ji} - \sum_{g\in\calG_i} p_{g} + p_{i}^{d} = 0, \quad \forall i \in \mathcal{B}\label{eq:powerbalance_1}\\
& -B_{ii} (e_i^2 + f_i^2) + \sum_{j \in \mathcal{N}_i^\text{fr}} q_{ij} + \sum_{j \in \mathcal{N}_i^\text{to}} q_{ji} - \sum_{g\in\calG_i} q_{g} + q_{i}^{d} = 0, \quad \forall i \in \mathcal{B} \label{eq:powerbalance_2}\\
& p_{i j}=-G_{i j}\left(e_{i}^{2}+f_{i}^{2}-e_{i} e_{j}-f_{i} f_{j}\right)-B_{i j}\left(e_{i} f_{j}-e_{j} f_{i}\right), \quad \forall (i,j) \in \mathcal{L}\label{eq:powerflow_1}\\
& p_{j i}=-G_{j i}\left(e_{j}^{2}+f_{j}^{2}-e_{j} e_{i}-f_{j} f_{i}\right)-B_{j i}\left(e_{j} f_{i}-e_{i} f_{j}\right), \quad \forall (i,j) \in \mathcal{L}\\
& q_{i j}=B_{i j}\left(e_{i}^{2}+f_{i}^{2}-e_{i} e_{j}-f_{i} f_{j}\right)-G_{i j}\left(e_{i} f_{j}-e_{j} f_{i}\right), \quad \forall (i,j) \in \mathcal{L}\\
& q_{j i}=B_{j i}\left(e_{j}^{2}+f_{j}^{2}-e_{j} e_{i}-f_{j} f_{i}\right)-G_{j i}\left(e_{j} f_{i}-e_{i} f_{j}\right), \quad \forall (i,j) \in \mathcal{L}\label{eq:powerflow_4}\\
& p_{i j}^{2}+q_{i j}^{2} \leq \bar{s}_{i j}^{2}, 
\quad p_{j i}^{2}+q_{j i}^{2} \leq \bar{s}_{i j}^{2}, \quad \forall(i, j) \in \mathcal{L},\label{eq:linelimits} \\
& \underline{v}_{i}^{2} \leq e_{i}^{2}+f_{i}^{2} \leq \bar{v}_{i}^{2}, \quad \forall i \in \mathcal{B}\label{eq:voltagelimit}\\
& \underline{p}_g \leq p_g \leq \overline{p}_g, 
\quad \underline{q}_g \leq q_g \leq \overline{q}_g, \quad \forall g \in \mathcal{G},\label{eq:genlimits}
\end{align}
\label{eq:acopf}
\end{subequations}

\vspace{-10pt}
The input to the optimization program is the power demands $x=\{p_i^d,q_i^d:i\in\mathcal{B}\}\in\mathbb{R}^{2|\mathcal{B}|}$. The output is the decision variables $p_g,q_g,e_i,f_i$, which are heavily coupled. When real power $p_g$ and voltage magnitude $v_i$ are given for all $g\in\calG,i\in\calB_\text{PV}$, the rest of the decision variables can be uniquely determined by the following system of power flow equations:
\vspace{-5pt}
\begin{align}
\label{eq:PF}
    \eqref{eq:powerbalance_1}-\eqref{eq:powerflow_4}, \quad
    e_{i}^{2}+f_{i}^{2} = v_{i}^{2}, \quad \forall i \in \calB_\text{PV}.
\end{align}


\vspace{-5pt}
We name $y_c=((p_g)_{g\in\calG_i},v_i)_{i\in\calB_\text{PV}}$ the control variables, as the specification of $y_c$ is sufficient for controlling the operation of the power system. We denote by $y_s$ the other decision variables of Eq.~\eqref{eq:acopf} and refer to them as state variables. Letting $\mathcal{Y} \subseteq \mathbb{R}^{|\mathcal{G}|+|\calB_\text{PV}|}$ and $\mathcal{S} \subseteq \mathbb{R}^{|\mathcal{G}|+2|\mathcal{B}|+4|\mathcal{L}|}$ denote the space of control and state variables, we define $PF:\mathcal{Y}\times \mathbb{R}^{2|\mathcal{B}|}\rightarrow\mathcal{S}$ as the mapping from control variables and input demands  to the state variables as the solution of Eq.~\eqref{eq:PF}.
Under this notation, we can rewrite the ACOPF objective in Eq.~\eqref{eq:acopf} as 
\begin{align}
\begin{aligned}
    \min_{y_c\in\mathcal{C}_c,y_s\in \mathcal{C}_s} \sum_{g\in\mathcal{G}} ( c_{2,g} p_g^2 + c_{1,g} p_g )
    \quad \text{subject to}\quad y_s \in PF(y_c, x)
\end{aligned}\label{eq:acopf_abstraction}
\end{align}
where $\mathcal{C}_c\subseteq\mathcal{Y}$ and $\mathcal{C}_s\subseteq\mathcal{S}$ represent the (convex) feasibility sets of control and state variables, respectively, as follows:
\begin{gather*}
    \mathcal{C}_c = \{((p_g)_{\calG_i},v_i)_{i\in\calB_\text{PV}}:\underline{p}_g \leq p_g \leq \overline{p}_g, \; \forall g\in\calG, \quad \underline{v}_i\leq v_i\leq\overline{v}_i, \; \forall i\in\calB_\text{PV}\},\\
    \hspace{-200pt}\mathcal{C}_s = \Big\{((q_i^g)_{i\in\calG_i},e_i,f_i,p_{ij},q_{ij},p_{ji},q_{ji})_{i\in\mathcal{B}}:\\
    \hspace{100pt}\underline{q}_g \leq q_g \leq \overline{q}_g, \; \forall g\in\calG, \; p_{ij}^2+q_{ij}^2\leq\overline{s}_{ij}^2,\; p_{ji}^2+q_{ji}^2\leq\overline{s}_{ij}^2,\,\forall (i,j)\in\calL\Big\}.
\end{gather*}
Note that $PF(y_c,x)$ many have a unique solution, multiple solutions, or no solution for arbitrary control variable $y_c$ and load $x$.


\section{A Constrained Machine Learning Approach to ACOPF}\label{sec:learn_acopf}

The (non-convex) QCQP in \eqref{eq:acopf_abstraction} can be stably solved by various existing algorithms including the interior-point method. However, as the size of power system scales up, the amount of computation required by an optimization solver becomes enormous. Considering the fact that the ACOPF problem needs to be repeatedly solved in real-life power systems as the power demand constantly changes, we are motivated to investigate alternative approaches that trade-off slight sub-optimality of the solution for computational efficiency.

\subsection{ML Approach for Control Prediction}

Given a fixed power system, our aim is to design a data-driven learning-based method that leverages samples of paired input and output of the ACOPF problem solved for the specific system under varying loads.
Suppose we have $N$ sample pairs $\{(x^n,y^n):n=1,\dots,N\}$ where $x^n$ denotes the power demands from the $n_{\text{th}}$ sample and $y^n$ denotes the optimal solution of \eqref{eq:acopf} under input $x^n$, which we can split into control and state variables $y^n=(y_c^n,y_s^n)$. A straightforward supervised learning framework for predicting the control variables $y_c^n$ from $x^n$ looks for a mapping $g:\mathbb{R}^{2|\mathcal{B}|}\rightarrow\mathcal{Y}$ that minimizes the data mismatch in the following way
\begin{align}
    \min_{g}\quad\sum_{n=1}^{N}\|g(x^n)-y_c^n\|^2.
\end{align}
We parameterize the mapping $g$ by a feed-forward neural network.

Despite its simplicity, this objective does not enforce constraint satisfaction, especially on the state variable. Specifically, $PF\big(g(x^n),x^n\big)$, the state variable resulting from the predicted control, may not lie within the constraint set $\mathcal{C}_s$. Such infeasible solutions 
require post-processing before they can be deployed; even after post-processing, the solutions are not guaranteed to be feasible and may degrade in generation cost.

To explicitly enforce the constraints to be satisfied, our work proposes a more sophisticated loss function that penalizes both data mismatch and constraint violation. Our initial proposal is to solve
\begin{align}
    \min_{g}\quad&\sum_{n=1}^{N}\Big(\|g(x^n)-y_c^n\|^2+ w \cdot r_{\mathcal{C}_s}\big(PF(g(x^n),x^n)\big)\Big)\label{eq:obj_with_PF}
\end{align}
where $w$ is a given weight value (e.g. $0.1$ and $1.0$ for our numerical experiment), and $r_{\mathcal{C}_s}$ is the following penalty associated with the violation of constraint set $\mathcal{C}_s$
\begin{align*}
    r_{\mathcal{C}_s} \hspace{-1pt}\Big((q_g)_{g\in\calG},(e_i,f_i)_{i\in\mathcal{B}},(p_{ij},q_{ij})_{ij\in\mathcal{L}},(p_{ji},q_{ji})_{ij\in\mathcal{L}}\Big)\hspace{-2pt}=\hspace{-2pt} \sum_{g\in\mathcal{G}}\hspace{-2pt}\left(\max\{0, \underline{q}_g\hspace{-2pt}-\hspace{-2pt}q_g\} + \max\{0,q_g-\overline{q}_g\}\right)&\\
    +\sum_{(i,j)\in\mathcal{L}}\left(\max\{0,p_{ij}^2+q_{ij}^2-\overline{s}_{ij}^2\} + \max\{0,p_{ji}^2+q_{ji}^2-\overline{s}_{ij}^2\}\right).&
\end{align*}

\subsection{Infeasible Power Flow System}

Training under the loss function in \eqref{eq:obj_with_PF}, however, may be unstable. The challenge comes from the fact that we may not always find a feasible solution $PF(g(x^n),x^n)$. This happens either because $PF(g(x^n),x^n)=\emptyset$ under control prediction $g(x^n)$ and/or load sample $x^n$, or because a numerical method (e.g. Newton-Raphson) cannot find a solution due to the lack of convergence guarantee. Training stagnates when this issue arises and the entire learning pipeline can be broken. To address the issue and stabilize training, we introduce slack variables $\sigma\in\mathbb{E}^{2|\mathcal{B}|}$ that captures the 
the minimum gap in power demand satisfaction for the power flow equations to have a solution.
Under the control prediction $g(x^n)$ and demand load $x^n$, we find approximate state variables by solving the following optimization program (which can be written as a non-convex QCQP after a simple reformulation of the objective):
\begin{align}
\label{eq:relaxed_PF}
    (\widehat{y}_s^n,\widehat{\sigma}^n) \in\operatorname{argmin}_{y_s\in\mathcal{S}, \sigma\in\mathbb{R}^{2|\mathcal{B}|}} \|\sigma\|_1
    \quad \text{subject to}\quad y_s \in PF(g(x^n), x^n+\sigma)
\end{align}
To stably learn a constrained ACOPF solution, we ultimately solve the bi-level optimization problem 
\begin{align}
    \min_{g}\quad&L(g)\triangleq\sum_{n=1}^{N}\left(\|g(x^n)-y_c^n\|^2+w\cdot r_{\mathcal{C}_s}\big(\widehat{y}_s^n\big)\right)\label{eq:bi_level}
\end{align}
where $\widehat{y}_s^n$ is defined in the lower-level optimization problem \eqref{eq:relaxed_PF}.
The loss function is tailored to the learning of ACOPF solutions leveraging the structure of the problem and is novel in the literature to the best of our knowledge. From a computational perspective, this loss, nevertheless, introduces significant challenges. Since $\widehat{y}_s^n$ is a function of $g(x^n)$ only defined implicitly through \eqref{eq:relaxed_PF}, it is unclear how the gradient of $r_{\mathcal{C}_s}\big(\widehat{y}_s^n\big)$ can be computed with respect to $g$, or even more fundamentally, whether $r_{\mathcal{C}_s}\big(\widehat{y}_s^n\big)$ is differentiable. In the next section, we provide an affirmative answer to the question and present a systematic method for deriving the (sub)gradient of $L$ with respect to $g(x^n)$ by adapting techniques from differentiable programming. Being able to evaluate this (sub)gradient means that we can compute the $\nabla_g L(g)$ through the chain rule, which allows us to optimize $L(g)$ using first-order algorithms.

\section{Differentiable QCQP}
\label{sec:differentiable-qcqp}

In this section, we show that under the second-order sufficient condition
any QCQP of the form
\begin{subequations}
\label{eq:qcqp}
\begin{align}
    z^*=\min_{z\in\mathbb{R}^k} \quad & \frac{1}{2} z^\top P_0 z + q_0^\top z \label{eq:qcqp:obj} \\
    \text{s.t.} \quad
    & \frac{1}{2} z^\top P_i z + q_i^\top z + r_i \leq 0 \quad \forall i=1,\dots,m_I, \\
    & \frac{1}{2} z^\top D_i z + h_i^\top z + g_i = 0 \quad \forall i=1,\dots,m_E
\end{align}
\end{subequations}
defines a \textit{differentiable} mapping from the parameters $P_0\in\mathbb{R}^{k\times k},q_0\in\mathbb{R}^k,\{P_i\in\mathbb{R}^{k\times k}\},\{q_i\in\mathbb{R}^k:i=1,\cdots,m_I\},\{r_i:\in\mathbb{R}:i=1,\cdots,m_I\}, \{D_i\in\mathbb{R}^{k\times k}\},\{ h_i\in\mathbb{R}^k:i=1,\cdots,m_E\}, \{g_i:\in\mathbb{R}:i=1,\cdots,m_E\}$ to the optimal solution $z^*$. We derive the gradient of $\ell(z^*)$ with respect to these parameters, where $\ell:\mathbb{R}^k\rightarrow\mathbb{R}$ can be any downstream loss function on $z^*$. When the second-order sufficient condition does not hold, we derive the subgradients within the subdifferentials $\partial_{P_0} \ell$, $\partial_{q_0} \ell$, etc, which allows subgradient descent/ascent to be performed on the downstream loss function. We note that \eqref{eq:qcqp} is a more general problem that covers \eqref{eq:relaxed_PF} as a special case.

Inspired by the literature on differentiable quadratic programming \cite{amos2017optnet} and convex programming \cite{agrawal2019differentiable}, we exploit a key structure to drive the innovation -- the KKT equations of \eqref{eq:qcqp} are preserved at the optimal solution $(z^*,\nu^*,\lambda^*)$ under differential changes in the parameters, where $\nu^*\in\mathbb{R}^{m_I}$ and $\lambda^*\in\mathbb{R}^{m_E}$ are the optimal dual variables associated the inequality and equality constraints, respectively. More specifically, under differential changes $\{d P_i\},\{dq_i\},\{dr_i\},\{dD_i\},\{dh_i\},\{dg_i\}$, we can find how $dz^*$ will change accordingly by solving a linear system of equations of the form
\begin{align}
    M\left[(dz^*)^{\top}, (d\nu^*)^{\top}, (d\lambda^*)^{\top}\right]^{\top} = b\Big(\{d P_i\},\{dq_i\},\{dr_i\},\{dD_i\},\{dh_i\},\{dg_i\}\Big),\label{eq:KKT_differential_system}
\end{align}
where we define in the appendix the matrix $M\in\mathbb{R}^{(k+m_I+m_E)\times(k+m_I+m_E)}$, which only depends on the parameters and optimal solution of \eqref{eq:qcqp}, and the vector $b(\{d P_i\},\{dq_i\},\{dr_i\},\{dD_i\},\{dh_i\},\{dg_i\})\in\mathbb{R}^{k+m_I+m_E}$, which is a function of the differentials of parameters. Setting the differentials of the parameters to appropriate identity matrices/tensors and solving this system of equations give the partial derivatives $\frac{\partial z^*}{\partial P_i}$, $\frac{\partial z^*}{\partial q_i}$, etc.,  by definition.

When $z^*$ is used to compute a differentiable downstream loss $\ell(z^*)$, we design a computationally efficient method for propagating the gradient of the loss through the QCQP to all parameters.
We state the main results below and defer the detailed derivation to the appendix.

\begin{thm}\label{thm:main}
    Suppose that strict complementary slackness, linear constraint qualification and second-order sufficient conditions hold at $(x^*,\nu^*,\lambda^*)$. Then, the matrix $M$ is invertible and $\ell$ is a differentiable function of the QCQP parameters. Given $\frac{\partial \ell}{\partial z^*}$, we have 
    \begin{align}
    \begin{aligned}
    &\nabla_{P_0} \ell = z^* d_z^{\top},\quad \nabla_{q_0} \ell = d_z,\quad\nabla_{P_i} \ell = \nu_i^* z^* d_z^{\top} + \frac{1}{2}\nu_i^* z^* (z^*)^{\top} d_{\nu_i},\quad\nabla_{q_i} \ell = \nu_i^* d_z + z^* d_{\nu_i} \\
    &\nabla_{r_i} \ell= d_{\nu_i}, \quad\nabla_{D_i} \ell = \lambda_i^* z^* d_z^{\top} + \frac{1}{2}d_{\lambda_i} z^* (z^*)^{\top},\quad\nabla_{h_i} \ell = \lambda_i^* d_z + z^* d_{\lambda_i},\quad\nabla_{g_i} \ell= d_{\lambda_i},
    \end{aligned}
    \label{thm:main:eq1}
    \end{align}
    where $d_z\in\mathbb{R}^{k}, d_{\nu}\in\mathbb{R}^{m_I}, d_{\lambda}\in\mathbb{R}^{m_E}$ are the solutions to
    \begin{align*}
    \left[ d_z^{\top},\,\,
        d_{\nu} ^{\top},\,\,
        d_{\lambda}^{\top}\right]^{\top}=-M^{-\top}\Big[ \big(\frac{\partial \ell}{\partial z^*}\big)^{\top},\,\,\cdots\,\,0\,\,\cdots, \,\,\cdots \,\,0 \,\,\cdots\Big]^{\top}.
    \end{align*}
\end{thm}

In the theorem, we state a sufficient condition on the differentiability of the QCQP and provide a systematic way of computing the gradients of the downstream loss on $z^*$ with respect to the QCQP parameters. We note that when the assumptions of Theorem~\ref{thm:main} do not hold, $M$ may not be invertible, and the QCQP in general is not differentiable. In that case, the expressions in \eqref{thm:main:eq1} are the subgradients where $d_z,d_{\nu}, d_{\lambda}$ are any solutions of the the under-determined system
\begin{align}
    -M^{\top}\left[ d_z^{\top},\,\,d_{\nu} ^{\top},\,\,d_{\lambda}^{\top}\right]^{\top}=\Big[ \big(\frac{\partial \ell}{\partial z^*}\big)^{\top},\,\,\cdots\,\,0\,\,\cdots, \,\,\cdots \,\,0 \,\,\cdots\Big]^{\top}.\label{eq:singular_M}
\end{align}
In our work, we take the solution with the smallest $\ell_2$ norm.

\section{Numerical Experiments}
\label{sec:experiments}

In this section, we demonstrate the performance of the proposed approach for ACOPF. 
We do not assume that power flow has a solution. If no solution exists for the power flow, we find a solution that minimizes the violation of the power balance constraints, as proposed in \eqref{eq:relaxed_PF}.
We highlight that our approach is capable of training the neural network model even with control prediction that can cause no solution of the power flow system (i.e., $PF(g(x^n),x^n) = \emptyset$). Note that Newton-Rhapson method fails to converge for such cases.

Our numerical experiments aim to demonstrate that the model can be successfully trained even with a number of training and testing epochs with the samples of infeasible power flow system (i.e., no solution and thus positive slack values).



\subsection{Experiment Settings}

We generate the training and testing data sets by solving the problem with perturbation of the active and reactive loads by random numbers uniformly generated from $(-1,1)$.
We consider three IEEE test instances, each of which has 30, 118, and 300 buses, taken from PGLib-OPF v21.07.
Because some problems can be infeasible with the random perturbations, we use the problem formulation that penalizes the violation of the power balance constraints.
Each sample $s$ of the data set consists of $p_s^d,q_s^d,p_g^*,q_g^*,e_i^*,f_j^*$, where $p_g^*,q_g^*,e_i^*,f_j^*$ are local optimal solution obtained by Ipopt v3.14.12. The training data set has 10,000 samples, and the testing data set has 2,500 samples.

We have implemented the proposed approach by using PyTorch v2.0.1, where the optimization problem is modeled with Pyomo v6.7.0 and solved by Ipopt v3.14.12. The trainings were parallelized with Horovod v0.28.1.
All the experiments were run on a Linux workstation with 144 CPUs of Intel(R) Xeon(R) Gold 6140 CPU @ 2.30GHz.

We use the feed-forward neural network with two layers and the sigmoid function as output. We use the same architectures proposed in \cite{pan2020deepopf}. Detailed experiment settings are given in Table~\ref{tab:settings}. Recall that the input and output dimensions are $|\calB_\text{PV}|+|\calG|$ and $2|\calB|$, respectively. All the models were trained by Adam optimizer with the learning rates given in the table.

\begin{table}[ht]
    \centering
    \begin{tabular}{cccccccc}
        IEEE Test System & $\calB$ & $\calB_\text{PV}$ & $\calG$ & $\calL$ & \# parameter per layer & $w$ & learning rate\\
        \hline
        30-bus      & 30  & 5  & 6  & 41  & 64, 32    & 1.0 & $10^{-4}$ \\
        118-bus     & 118 & 53 & 54 & 231 & 256, 128  & 1.0 & $10^{-4}$\\
        300-bus     & 300 & 68 & 69 & 411 & 1024, 512 & 0.1 & $10^{-6}$\\
    \end{tabular}
    \caption{Experiment setting for each IEEE test system}
    \vspace{-15pt}
    \label{tab:settings}
\end{table}

\subsection{Small Test Cases}

\begin{figure}[ht]
\centering
    \subfigure[30-bus network]{
        \includegraphics[width=.3\textwidth]{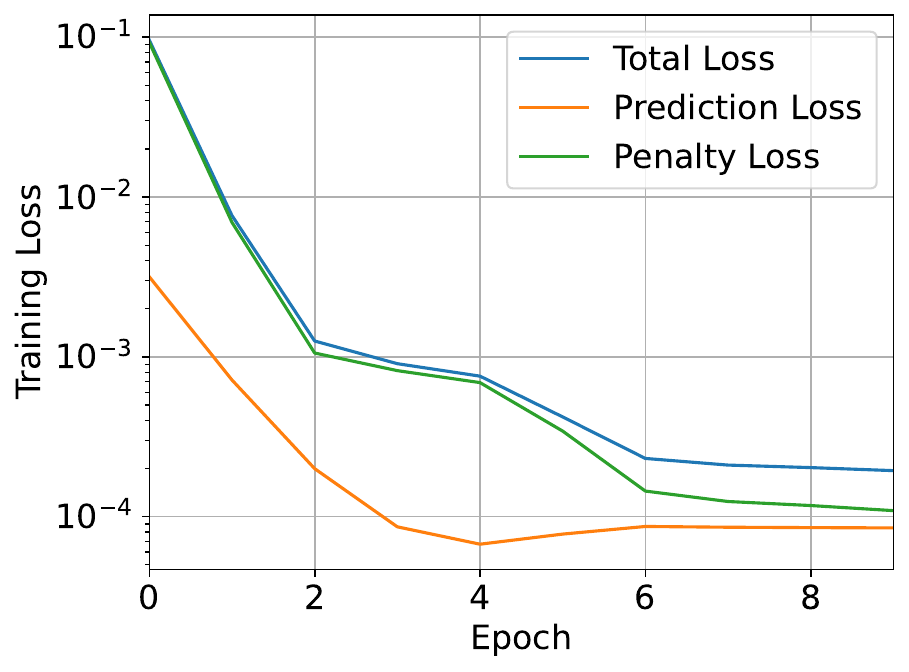}
    }
    \subfigure[118-bus network]{
        \includegraphics[width=.3\textwidth]{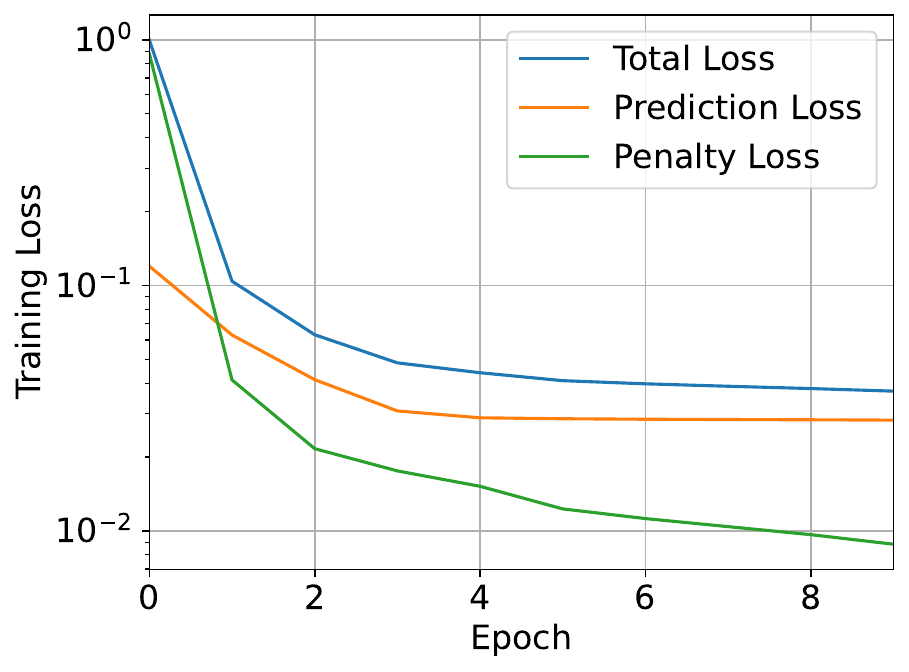}
    }
    \subfigure[Testing accuracy]{
        \includegraphics[width=.3\textwidth]{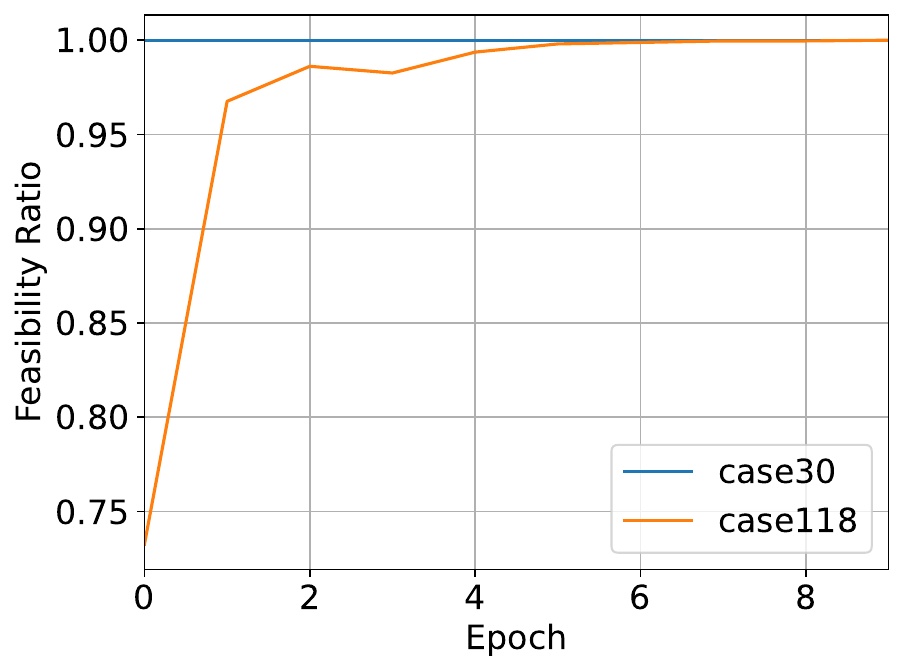}
    }
    \vspace{-5pt}
    \caption{Performances on small grid networks (IEEE 30- and 118-bus systems)}
    \vspace{-12pt}
    \label{fig:small}
\end{figure}

Figure~\ref{fig:small} presents the training and testing performances of the proposed approach for the small test networks: IEEE 30- and 118-bus network systems.
The prediction loss and penalty loss values in the figure measure the first and second terms of \eqref{eq:bi_level}. The total loss measures $L(g)$.
The 100\% test accuracy (i.e., feasibility ratio) has been obtained within 10 epochs for both small systems; that is, the control predictions $g(x^n)$ made from the trained models result in feasible state solutions with respect to the power flow equations \eqref{eq:PF} for all 2,500 test samples.
Our results are consistent to that presented in \cite{pan2020deepopf}, where the power flow equations of \eqref{eq:PF} are solved as compared to the optimization \eqref{eq:relaxed_PF} in our approach.

\subsection{Large Test Case with Infeasible Power Flow}

\begin{figure}[ht]
\centering
    \subfigure[Training loss]{
        \includegraphics[width=.3\textwidth]{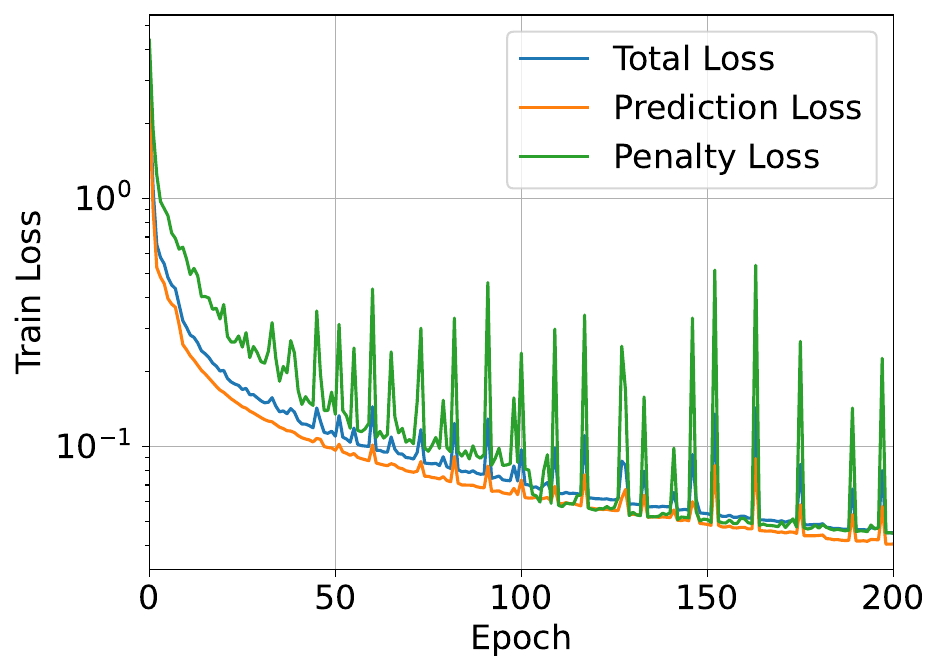}
    }
    \subfigure[Infeasible power flow]{
        \includegraphics[width=.3\textwidth]{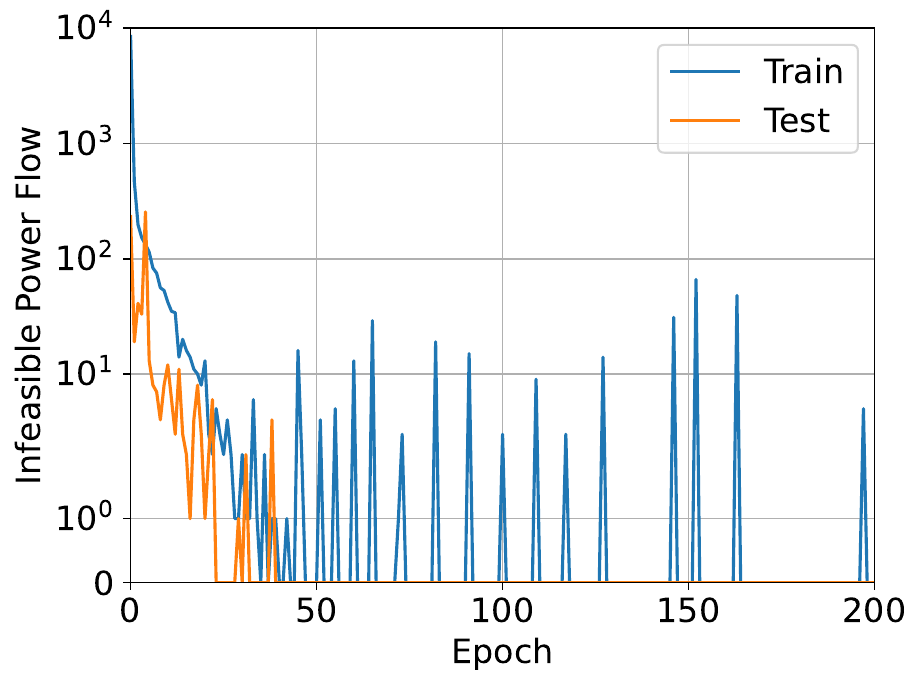}
    }
    \subfigure[Testing accuracy]{
        \includegraphics[width=.3\textwidth]{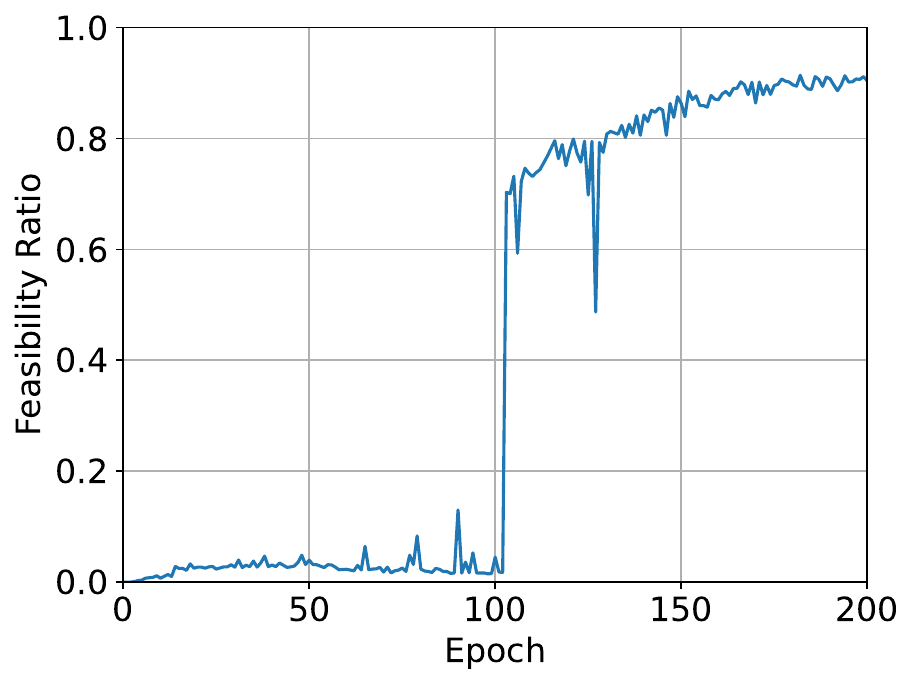}
    }
    \vspace{-5pt}
    \caption{Training and testing performances on the IEEE 300-bus network}
    \label{fig:case300}
\end{figure}

Next, we demonstrate that the neural network model can be successfully trained even when some control prediction leads to infeasible power flow equations.
In Figure~\ref{fig:case300} we report the training loss, the number of infeasible power flow solves, and the accuracy of the model over the training epochs for the IEEE 300-bus network data. Recall that \eqref{eq:PF} has no solution if the values of the slack variable $\sigma$ in \eqref{eq:relaxed_PF} are positive.
The model trained on IEEE 300-bus network system achieves a test accuracy of 92\% after 200 epochs.
In a number of training and testing epochs (Figure \ref{fig:case300}b), we observe that the neural network predicts control variables that result in the infeasible state variables. 
While reducing over the epochs as the neural network learns more accurate controls, a positive number of infeasible power flow cases still appear in many of the later epochs. We observe that the training epochs with such infeasible cases experience the spikes in the penalty loss value.
We highlight that our approach with the QCQP activation function (with slack variables) allows the neural network to still learn and improve when such infeasible cases arise, while existing approaches relying on the power flow equation solver (e.g. \cite{pan2020deepopf}) do not.
\vspace{-5pt}
\vspace{-15pt}


\section{Concluding Remarks}\label{sec:conclusion}

Our work proposes a learning-based framework for solving the notoriously challenging ACOPF problem with the aim of achieving computational efficiency and reliably generating high-quality feasible solutions. 
We identify that a common issue of the existing approaches in this domain \cite{pan2020deepopf, donti2021dc3} lies in the assumption that the power flow equation admits a solution for any neural-network-predicted control variables, which does not always hold. Training gets disrupted when the power flow solver fails to produce a feasible solution. We address this issue by modeling the power flow as a non-convex QCQP problem that minimizes the constraint violation. By leveraging and generalizing techniques from differentiable convex programming, we derive (sub)gradient of the state variables with respect to the control variables, which allows the loss function on state variables to be properly back-propagated. We show through numerically simulations that our proposed framework stably learns solutions with high feasibility rate and low generation in large systems with wide load variations, in which existing approaches fail to train.

\section*{Disclaimer}
This paper was prepared for informational purposes in part by
the Artificial Intelligence Research group of JP Morgan Chase \& Co and its affiliates (``JP Morgan''),
and is not a product of the Research Department of JP Morgan.
JP Morgan makes no representation and warranty whatsoever and disclaims all liability,
for the completeness, accuracy or reliability of the information contained herein.
This document is not intended as investment research or investment advice, or a recommendation,
offer or solicitation for the purchase or sale of any security, financial instrument, financial product or service,
or to be used in any way for evaluating the merits of participating in any transaction,
and shall not constitute a solicitation under any jurisdiction or to any person,
if such solicitation under such jurisdiction or to such person would be unlawful.

\bibliographystyle{plain} 
\bibliography{references}

\appendix

\clearpage
\appendix

\section{Appendix - Detailed Derivation of Results in Section~\ref{sec:differentiable-qcqp} \& Proof of Theorem~\ref{thm:main}}

In this section, we first show how we have derived the system of equations \eqref{eq:KKT_differential_system} that the differentials need to satisfy. We start by writing the Lagrangian of the optimization problem \eqref{eq:qcqp}.
\begin{align}
    L(z, \nu, \lambda)
    =\frac{1}{2} z^{\top} P_0 z+q_0^{\top} z
    +\sum_{i=1}^{m_I}\nu_i\left(\frac{1}{2} z^{\top} P_i z+q_i^{\top}z+r_i\right)
    +\sum_{i=1}^{m_E}\lambda_i\left(\frac{1}{2} z^{\top} D_i z+h_i^{\top}z+g_i\right),
\end{align}
Here $\nu\in\mathbb{R}^{m_I}$ and $\lambda\in\mathbb{R}^{m_E}$ are the Lagrangian multipliers of the inequality and equality constraints.

The KKT optimality conditions for stationarity, primary and dual feasibility, and complementary slackness are
\begin{subequations}
\begin{align}
    & P_0 z^*+q_0+\sum_{i=1}^{m}\nu_i^*\left(P_i z^*+q_i\right)+\sum_{j=1}^{n}\lambda_i^*\left(D_i z^* +h_i\right)=0 \\
    & \frac{1}{2} (z^*)^{\top} P_{i} x^*+q_{i}^{\top} z^*+r_{i} \leq 0 && \forall i=1, ..., m_I, \\
    &\frac{1}{2} (z^*)^{\top} D_{i} x^*+h_{i}^{\top} z^*+g_{i} = 0 && \forall j=1, ..., m_E, \\
    &\nu_i^*\left(\frac{1}{2} (z^*)^{\top} P_i z^*+q_i^{\top} z^*+r_i\right)=0 && \forall i=1, ..., m_I,\\
    &\nu_i^* \ge 0 && \forall i=1,\dots,m_I.
\end{align}
\end{subequations}

Taking the differentials of the equality equations, we get the following system of equations
\begin{subequations}
\label{eq:KKT}
\begin{align}
    0 &= dP_0 z^*+P_0 dz^*+dq_0+\sum_{i=1}^{m}\left(d\nu_i^* P_i z^*+\nu_i^* dP_i z^*+\nu_i^* P_i dz^*+d\nu_i^* q_i+\nu_i^* dq_i\right)\notag\\
    & \qquad +\sum_{j=1}^{n}\left(d\lambda_i^* D_i z^*+\lambda_i^* dD_i z^*+\lambda_i^* D_i dz^* + d\lambda_i^* h_i+\lambda_i^* dh_i\right), \\
    0 &= \frac{1}{2}(dz^*)^{\top}D_i z^* + \frac{1}{2}(z^*)^{\top}dD_i z^* + \frac{1}{2}(z^*)^{\top}D_i dz^*+dh_i^{\top}z^*+h_i^{\top}dz^*+dg_i,\,\, \forall i=1,\dots,m_E, \\
    0 &= \nu_i^*\left(\frac{1}{2 } (dz^*)^{\top} P_i z^*+\frac{1}{2}(z^*)^{\top} dP_i z^* + \frac{1}{2} (z^*)^{\top} P_i dz^* +dq_i^{\top} z^* + q_i^{\top} dz^* +dr_i\right) \notag \\
    & \qquad +d\nu_i^*\left(\frac{1}{2} (z^*)^{\top} P_{i} z^*+q_{i}^{\top} z^*+r_{i}\right),\,\, \forall i=1, \dots, m_I.
\end{align}
\end{subequations}

The system \eqref{eq:KKT} of equations can be written in the concise matrix form \eqref{eq:KKT_differential_system} where
\begin{equation*}
    \begin{aligned}
        M &= \left[\begin{array}{ccccc}
        P_0+\sum_{i=1}^{m}\nu_i^* P_i+\sum_{j=1}^{n}\lambda_i^* D_i & P_1 z^*+q_1 & \cdots & D_1 z^*+h_1 & \cdots \\
        \nu_1^*\left(\frac{1}{2}(z^*)^{\top}P_1^{\top}+\frac{1}{2}(z^*)^{\top} P_1 + q_1^{\top}\right) & \frac{1}{2}(z^*)^{\top} P_1 z^* + q_1^{\top} x^* + r_1 & \cdots & 0 & \cdots\\
        \vdots & \vdots &  & \vdots &  \\
        \frac{1}{2}(z^*)^{\top}D_1^{\top}+\frac{1}{2}(z^*)^{\top}D_1+h_1^{\top} & 0 & \cdots & 0 & \cdots \\
         \vdots & \vdots & & \vdots & 
        \end{array}\right].
    \end{aligned}
\end{equation*}
\begin{align}
    \text{and }&b\Big(\{d P_i\},\{dq_i\},\{dr_i\},\{dD_i\},\{dh_i\},\{dg_i\}\Big)\notag\\
    &\qquad=-\left[\begin{array}{c} dP_0 z^* +dq_0 + \sum_{i=1}^{m}(\nu_i^* dP_i z^*+\nu_i^* dq_i) + \sum_{i=j}^{n}(\lambda_i^* dD_i z^*+\lambda_i^* dh_i)\\
    \nu_1^*\left(\frac{1}{2}(z^*)^{\top} dP_1 z^* + dq_1^{\top} z^* + dr_1\right)\\
    \vdots\\ 
    \frac{1}{2}(z^*)^{\top}dD_1 z^*+dh_1^{\top} z^* + dg_1\\
    \vdots\end{array}\right].
\end{align}

As we have discussed in Section~\ref{sec:differentiable-qcqp}, we can find the partial derivative of $z^*$ with respect to the parameters by solving \eqref{eq:KKT_differential_system} with properly selected differentials in the parameters. For example, to compute $\frac{\partial z^*}{\partial r_1}$, we can set $d r_1$ to 1 and all other differentials to 0, solve the following system, and extract $dz^*$ from the solution
\begin{align*}
    M\left[(dz^*)^{\top}, (d\nu^*)^{\top}, (d\lambda^*)^{\top}\right]^{\top} = b\Big(\{\vct{0}_{k\times k}\},\{\vct{0}_k\},(1,0,\cdots,0),\{\vct{0}_{k\times k}\},\{\vct{0}_k\},\{0\}\Big).
\end{align*}
This procedure needs to be repeatedly applied to find the partial derivative of $z^*$ with respect to other parameters, which incurs huge computational costs. Fortunately, we show that when a downstream loss function $\ell(z^*)$ is computed on $z^*$, we can much more efficiently back-propagate the gradient through the QCQP.

When $M$ is invertible, by the chain rule
\begin{align*}
    \frac{\partial \ell}{\partial r_1} &= \left(\frac{\partial \ell}{\partial z^*}\right)^{\top}\frac{\partial z^*}{\partial r_1}\\
    &= \left[\left(\frac{\partial \ell}{\partial z^*}\right)^{\top},\vct{0}^{\top},\vct{0}^{\top}\right] M^{-1} b\big(\{\vct{0}_{k\times k}\},\{\vct{0}_k\},(1,0,\cdots,0),\{\vct{0}_{k\times k}\},\{\vct{0}_k\},\{0\}\big)\\
    &=b\big(\{\vct{0}_{k\times k}\},\{\vct{0}_k\},(1,0,\cdots,0),\{\vct{0}_{k\times k}\},\{\vct{0}_k\},\{0\}\big)^{\top}M^{-\top}\left[\left(\frac{\partial \ell}{\partial z^*}\right)^{\top},\vct{0}^{\top},\vct{0}^{\top}\right]^{\top},
\end{align*}
where the second equality follows from the fact that the loss function does not depend on the dual variables.

It is important to note that computing the gradient of $\ell$ with respect to any other parameter also involves calculating $M^{-\top}\left[\left(\frac{\partial \ell}{\partial z^*}\right)^{\top},\vct{0}^{\top},\vct{0}^{\top}\right]^{\top}$; we just need to hit the product on the left by a different $b$ vector. Having observed and exploited this fact, we can show that once we compute
\begin{align*}
    \left[ d_z^{\top},\,\,
        d_{\nu} ^{\top},\,\,
        d_{\lambda}^{\top}\right]^{\top}=-M^{-\top}\left[ \Big(\frac{\partial \ell}{\partial z^*}\Big)^{\top},\,\,\cdots\,\,0\,\,\cdots, \,\,\cdots \,\,0 \,\,\cdots\right]^{\top},
\end{align*}
the gradients will be those given in \eqref{thm:main:eq1}.

As we have explained in the paragraph after Theorem~\ref{thm:main}, when $M$ is not invertible, \eqref{thm:main:eq1} gives us the subgradients of the downstream loss when we solve \eqref{eq:singular_M} for $d_z,d_{\nu}, d_{\lambda}$.

Next, we prove that the matrix $M$ is invertible (non-singular) when strict complementary slackness, linear constraint qualification, and second-order sufficient conditions are satisfied.

Similar to~\cite{amos2017optnet}, we show that there is a unique solution to the following system of equations:
\begin{equation}
    \begin{bmatrix}
        Q(z^*,\nu^*,\lambda^*) & \nabla P(z^*) & \nabla D(z^*)\\
        Diag(\nu^*)\nabla P(z^*)^T & Diag(P(z^*)) & 0\\
        \nabla D(z^*)^T & 0 & 0,
    \end{bmatrix}
    \begin{bmatrix}
    z\\ \nu \\ \lambda
    \end{bmatrix}=
    \begin{bmatrix}
    a\\ b \\ c
    \end{bmatrix},
\end{equation}
where $Q(z,\nu,\lambda)$ and $Diag(\cdot)$ are the Hessian matrix of~\eqref{eq:qcqp} evaluated at $(z,\nu,\lambda)$ and a diagonal operator, respectively.
Note that $\nabla P(z^*):=\nabla P(z)|_{z=z^*}$ and the same for $\nabla D(z^*)$.

Let $\mathcal{A}(z^*)$ be the active set at $z^*$ such that the inequality constraints are satisfied as equalities, $\mathcal{A}(z^*):=\{i=1\,\dots,m_I \mid P_i(z^*)=0\}$.
Since the strict complementary slackness holds, we have $\nu^*_i > 0$ for all $i \in \mathcal{A}(z^*)$ and $\nu^*_i=0$ for all $i \notin \mathcal{A}(z^*)$.
Similar to~\cite{amos2017optnet}, for $i \notin \mathcal{A}(z^*)$ we set $\nu_i=b_i/P_i(z^*)$.
Then we need to show that the following system of equations has a unique solution:
\begin{equation}
    \begin{aligned}
        \begin{bmatrix}
            Q(z^*,\nu^*,\lambda^*) & \nabla P(z^*)_{\mathcal{A}(z*)} & \nabla H(z^*)\\
            Diag(\nu^*)\nabla P(z^*)_{\mathcal{A}(z^*)}^T & 0 & 0\\
            \nabla H(z^*)^T & 0 & 0,    
        \end{bmatrix}
        \begin{bmatrix}
        z \\ \nu_{\mathcal{A}(z^*)} \\ \lambda
        \end{bmatrix}
        =
        \begin{bmatrix}
        a - \nabla P(z^*)_{\mathcal{A}(z^*)} \\ b_{\mathcal{A}(z^*)} \\ c
        \end{bmatrix}\\
        (\Leftrightarrow)\;\begin{bmatrix}
            Q(z^*,\nu^*,\lambda^*) & \nabla P(z^*)_{\mathcal{A}(z^*)} & \nabla D(z^*)\\
            \nabla P(z^*)_{\mathcal{A}(z^*)}^T & 0 & 0\\
            \nabla D(z^*)^T & 0 & 0,    
        \end{bmatrix}
        \begin{bmatrix}
        z \\ \nu_{\mathcal{A}(z^*)} \\ \lambda
        \end{bmatrix}
        =
        \begin{bmatrix}
        a - \nabla P(z^*)_{\mathcal{A}(z^*)} \\ b_{\mathcal{A}(z^*)}/Diag(\nu^*)_{\mathcal{A}(z^*)} \\ c
        \end{bmatrix}
    \end{aligned}
\end{equation}
By linear constraint qualification, the matrix $K=\begin{bmatrix}\nabla P(z^*)_{\mathcal{A}(z^*)}^T \\ \nabla D(z^*)^T\end{bmatrix}$ has full row rank.
The second-order sufficient conditions ensure that the Hessian matrix $Q$ is positive-definite in the null space of $K$.
Therefore, the left-hand side matrix is non-singular\footnote{By multiplying $z$ to the first row of the matrix we have $z^TQz=0 \Rightarrow z=0$ and $\nu_{\mathcal{A}(z^*)}=\lambda=0$ follows from full rank assumption.}, thus having a unique solution.

\end{document}